\documentclass[12pt,leqno]{amsart}
\usepackage{amssymb}
\usepackage{amsmath,amssymb,color}
\numberwithin{equation}{section}

\newcommand{\sump}{\sum_{\ell=1}^p}

\newcommand{\la}{\lambda}
\newcommand{\va}{\varphi}
\newcommand{\ppp}{\partial}

\newcommand{\pppa}{\partial_t^{\alpha}}

\newcommand{\R}{\mathbb{R}}
\newcommand{\Q}{\mathbb{Q}}

\newcommand{\C}{\mathbb{C}} 
\newcommand{\N}{\mathbb{N}}

\newcommand{\ooo}{\overline}
\newcommand{\OOO}{\Omega}
\newcommand{\MLO}{E_{\alpha,1}}
\newcommand{\MLT}{E_{\alpha,2}}

\newcommand{\sumn}{\sum_{n=1}^{\infty}}
\newcommand{\sumk}{\sum_{k=1}^{d_n}}
\newcommand{\sumij}{\sum_{i,j=1}^d}

\newcommand{\hhalf}{\frac{1}{2}}
\newcommand{\DDD}{\mathcal{D}}

\allowdisplaybreaks

\title
[]
{
Decay rates and initial values for time-fractional diffusion-wave equations
}

\author{
$^{1,2,3}$ M.~Yamamoto }

\thanks{
$^1$ Graduate School of Mathematical Sciences, The University
of Tokyo, Komaba, Meguro, Tokyo 153-8914, Japan \\
$^2$ Honorary Member of Academy of Romanian Scientists, 
Splaiul Independentei Street, no 54,
050094 Bucharest Romania \\
$^3$ Peoples' Friendship University of Russia 
(RUDN University) 6 Miklukho-Maklaya St, Moscow, 117198, Russian Federation
e-mail: {\tt myama@ms.u-tokyo.ac.jp}
}

\date{}
\begin{document}
\maketitle

\baselineskip 18pt

\begin{abstract}
We consider a solution $u(\cdot,t)$ to an
initial boundary value problem for time-fractional diffusion-wave 
equation with the order $\alpha \in (0,2) \setminus \{ 1\}$ where 
$t$ is a time variable.
We first prove that a suitable norm of $u(\cdot,t)$ is bounded by   
$\frac{1}{t^{\alpha}}$ for $0<\alpha<1$ 
and $\frac{1}{t^{\alpha-1}}$ for $1<\alpha<2$ for all large $t>0$.
Moreover we characterize initial values in the cases where the decay rates
are faster than the above critical exponents.  Differently from the 
classical diffusion equation $\alpha=1$, the decay rate can give some 
local characterization of initial values.
The proof is based on the eigenfunction expansions of solutions and the 
asymptotic expansions of the Mittag-Leffler functions for large
time.
\end{abstract}

\section{Introduction}

Let $\OOO \subset \R^d$ be a bounded domain with smooth boundary
$\ppp\OOO$ and let $\nu(x) := (\nu_1(x), ...., \nu_d(x))$ be the 
unit outward normal vector to $\ppp\OOO$ at $x$.
We assume that 
$$
0< \alpha < 2, \quad \alpha \ne 1.
$$
By $\pppa$ we denote the Caputo derivative:
$$
\pppa g(t) = \frac{1}{\Gamma(n-\alpha)}\int^t_0
(t-s)^{n-\alpha-1}\frac{d^n}{ds^n}g(s) ds
$$
for $\alpha \not\in \N$ satisfying $n-1<\alpha<n$ with 
$n\in \N$ (e.g., Podlubny \cite{Po}).  For $\alpha=1$, 
we write  
$\ppp_tg(t) = \frac{dg}{dt}$ and $\ppp_tg(x,t) = \frac{\ppp g}{\ppp t}
(x,t)$.

We consider an initial boundary value problem for a time-fractional
diffusion-wave equation: 
$$
\left\{ \begin{array}{rl}
& \pppa u(x,t) = -Au(x,t), \quad x\in \OOO, \, 0<t<T, \\
& u\vert_{\ppp\OOO\times (0,T)} = 0, \\
& u(x,0) = a(x), \quad x \in \OOO \quad \mbox{if $0<\alpha \le 1$}, \\
& u(x,0) = a(x), \quad \ppp_tu(x,0) = b(x), \quad x \in \OOO 
\quad \mbox{if $1<\alpha < 2$}.
\end{array}\right.
                                       \eqno{(1.1)}
$$
Throughout this article, we set  
$$
(-Av)(x) = \sumij \ppp_i(a_{ij}(x)\ppp_jv(x)) + c(x)v(x), \quad 
x\in \OOO,
$$
where $a_{ij} = a_{ji}$, $1\le i,j\le n$ and $c$ are sufficiently smooth 
on $\ooo{\OOO}$, and $c(x) \le 0$ for $x \in \ooo{\OOO}$, and we
assume that there exists a constant $\sigma>0$ such that 
$$
\sumij a_{ij}(x)\zeta_i\zeta_j \ge \sigma \sum_{i=1}^d \zeta_i^2 
\quad \mbox{for all $x \in \ooo{\OOO}$ and $\zeta_1, ..., \zeta_d \in \R$}.   
$$

For $\alpha \in (0,2) \setminus \{1\}$, the first equation in (1.1) is 
called a fractional diffusion-wave equation, which models anomalous 
diffsion in heterogeneous media. 
As for physical backgrounds, we are restricted to a few references: 
Metzler and Klafter \cite{MK}, Roman and Alemany \cite{RA}, and 
one can consult Chapter 10 in \cite{Po}.
\\

The properties such as asymptotic behavior as $t \to \infty$ of 
solution $u$ to (1.1) are proved to depend on the fractional order $\alpha$ 
of the derivative.  Moreover decay rates can characterize the initial values
which is very different from the case $\alpha=1$.
The main purpose of this article is to study these topics.

Throughout this article, $L^2(\OOO)$, $H^{\mu}(\OOO)$ denote the usual 
Lebesgue space and Sobolev spaces (e.g., Adams \cite{Ad}), 
and by $\Vert \cdot\Vert$ and 
$(\cdot, \cdot)$ we denote the norm and the scalar product in 
$L^2(\OOO)$ respectively.  When we specify the norm in a Hilbert space
$Y$, we write $\Vert \cdot\Vert_Y$.
All the functions under consideration are assumed 
to be real-valued.
 
We define the domain $\mathcal{D}(A)$ of $A$ by $H^2(\OOO) \cap H^1_0(\OOO)$.
Then the operator $A$ in $L^2(\OOO)$ has positive eigenvalues 
with finite multiplicities.  We denote the set of all the eigenvalues by 
$$
0 < \la_1 < \la_2 \cdots \longrightarrow \infty.
$$
We set $\mbox{Ker}\, (A-\la_n):= \{ v\in \mathcal{D}(A);\, Av = \la_nv\}$
and $d_n:= \mbox{dim}\, \mbox{Ker}\, (A-\la_n)$. 
We denote an orthonormal basis 
of Ker $(A-\la_n)$ by $\{\va_{nk}\}_{1\le k \le d_n}$.

Then we define a fractional power $A^{\gamma}$ with 
$\gamma \in \R$ (e.g., Pazy \cite{Pa}), and we see
$$
\DDD(A^{\gamma}) = \left\{ a\in L^2(\OOO);\, 
\sumn \sumk \la_n^{2\gamma} (a,\va_{nk})^2 < \infty\right\} \quad 
\mbox{if $\gamma > 0$}
$$
and $\DDD(A^{\gamma}) \supset L^2(\OOO)$ if $\gamma \le 0$,
$$
A^{\gamma}a = \sumn \la_n^{\gamma} \sumk (a,\va_{nk})\va_{nk}, \quad
\Vert A^{\gamma}a\Vert = \left( \sumn \sumk \la_n^{2\gamma} (a,\va_{nk})^2
\right)^{\hhalf}, \quad a\in \DDD(A^{\gamma}).      \eqno{(1.2)} 
$$
In particular,
$$
A^{-1}a = \sumn \sumk \frac{1}{\la_n} (a,\va_{nk})\va_{nk},
\quad \Vert A^{-1}a \Vert 
= \left( \sumn \sumk \frac{1}{\la_n^2} (a,\va_{nk})^2\right)^{\hhalf}.
                                          \eqno{(1.3)}
$$
Moreover it is known that 
$$
\DDD(A^{\gamma}) \subset H^{2\gamma}(\OOO) \quad \mbox{for $\gamma \ge 0$}.
$$
The well-posedness for (1.1) is studied for example 
in Gorenflo, Luchko and Yamamoto \cite{GLY}, 
Kubica, Ryszewska and Yamamoto \cite{KRY}, Sakamoto and Yamamoto \cite{SY}.
As for the asymptotic behavior, we know 
$$
\Vert u(\cdot,t)\Vert \le \frac{C}{t^{\alpha}}\Vert a\Vert, \quad t>0
                                    \eqno{(1.4)}
$$
(e.g., \cite{KRY}, \cite{SY}, Vergara and Zacher \cite{VZ}).  
The article \cite{VZ} first established (1.4) for $t$-dependent operator $A$.
Moreover by the eigenfunction expansion of $u(x,t)$ (e.g., \cite{SY}), 
one can prove
$$
\Vert u(\cdot,t)\Vert \le \frac{C}{t^{\alpha}}\Vert a\Vert
+ \frac{C}{t^{\alpha-1}}\Vert b\Vert, \quad t>0               \eqno{(1.5)}
$$
for $1<\alpha<2$.

First we improve (1.4) and (1.5) with stronger norm of $u$.
\\
{\bf Theorem 1.}\\
{\it
Let $t_0>0$ be arbitrarily fixed.
There exists a constant $C>0$ depending on $t_0$ such that 
$$
\Vert u(\cdot,t)\Vert_{H^2(\OOO)} \le 
\left\{\begin{array}{rl}
& \frac{C}{t^{\alpha}}\Vert a\Vert \quad \mbox{if $0<\alpha<1$},\\
& \frac{C}{t^{\alpha}}\Vert a\Vert
+ \frac{C}{t^{\alpha-1}}\Vert b\Vert \quad \mbox{if $1<\alpha<2$}
\end{array}\right.
$$
for $t \ge t_0$.
}
\\

This theorem means that the Sobolev regularity of initial values 
is improved by $2$ after any time $t>0$ passes.

The fractional diffusion-wave equation (1.1) models slow diffusion, which  
the decay estimates (1.4) and (1.5) describe.  For $\alpha=1$, 
by the eigenfunction exansion of $u$, we can readily prove that  
$\Vert u(\cdot,t)\Vert \le e^{-\la_1t}\Vert a\Vert$.
Needless to say, Theorem 1 does not reject the exponential decay
$e^{-\la_1 t}$, but as this article shows, the decay rates in the theorem 
are the best possible in a sense. 

For further statements, we introduce a bounded linear operator 
$F: \DDD(A^{\gamma}) \, \longrightarrow\, Y$, where $\gamma > 0$ and
$Y$ is a Hilbert space with the norm $\Vert\cdot \Vert_Y$.
We interpret that $F$ is an observation mapping, and we 
consider the following four kinds of $F$.
\\
{\bf Case 1.}
\\
Let $\omega \subset \OOO$ be a subdomain.  Let
$$
F_1(v) = v\vert_{\omega}, \quad \DDD(F_1) = L^2(\OOO), \quad 
Y = L^2(\omega).                          \eqno{(1.6)}
$$
Then $F_1 : L^2(\OOO) \longrightarrow L^2(\omega)$ is bounded.
\\
{\bf Case 2.}
\\
Let $\Gamma \subset \ppp\OOO$ be a subboundary.  Let 
$$
F_2(v) = \ppp_{\nu_A}v\vert_{\Gamma}, \quad \DDD(F_2) = H^2(\OOO),
\quad Y = L^2(\Gamma).                           \eqno{(1.7)}
$$
Here we set 
$$
\ppp_{\nu_A}v:= \sumij a_{ij}(x)(\ppp_iv)(x)\nu_j(x).
$$
The trace theorem (e.g., Adams \cite{Ad}) implies that 
$F_2: H^2(\OOO) \longrightarrow L^2(\Gamma)$ is bounded.
\\
{\bf Case 3.}
\\
Let $x^1, ..., x^M \in \OOO$ be fixed and let $\gamma > \frac{d}{4}$, where 
$d$ is the spatial dimensions.
We consider
$$
F_3(v) = (v(x^1), ..., v(x^M)), \quad 
\DDD(F_3) = \DDD(A^{\gamma}), \quad Y = \R^M.       \eqno{(1.8)}
$$
Then the Sobolev embedding implies that $\DDD(F_3) \subset C(\ooo{\OOO})$,
and so $F_3: \DDD(A^{\gamma}) \longrightarrow \R^M$ is bounded.
We interpret that $F_3$ are pointwise data.
\\
{\bf Case 4.}
\\
Let $\rho_1, ..., \rho_M \in L^2(\OOO)$ be given and let $Y=\R^M$.
Let 
$$
F_4(v) = \left( \int_{\OOO} \rho_k(x)v(x) dx \right)_{1\le k \le M},
\quad \DDD(F_4) = L^2(\OOO), \quad Y = \R^M.        \eqno{(1.9)}
$$
Then $F_4: L^2(\OOO) \longrightarrow \R^M$ is 
bounded and corresponds to distributed data with 
weight functions $\rho_k$ whose supports concentrate around  
some points in $\OOO$.
\\

Now we state
\\
{\bf Theorem 2.}\\
{\it
In (1.1) we assume that $a, b \in L^2(\OOO)$ for 
$F_1, F_2, F_4$ and 
$a, b \in \DDD(A^{\gamma_0})$ with $\gamma_0 = 0$ if $\frac{d}{4} < 1$ and
$\gamma_0 > \frac{d}{4} - 1$ if $\frac{d}{4} \ge 1$ for $F_3$. 
Let $u = u(x,t)$ satisfy (1.1).  For $j=3,4$, let $F_j$ satisfy
$F_j\vert_{\mbox{Ker}\, (\la_n-A)}$ is injective for all $n\in \N$.

Firthermore we assume that for $j=1,2,3,4$, 
there exist sequences $\tau_n$, $n\in \N$ and 
$C_n>0$, $n\in \N$ which may depend on $u$, such that 
$$
\tau_n > 0, \quad \lim_{n\to\infty} \tau_n = \infty     \eqno{(1.10)}
$$
and
$$
\Vert F_j(u(\cdot,t))\Vert_Y \le \frac{C_n}{t^{\tau_n}}
\quad \mbox{as $t \to \infty$ for all $n\in \N$.}           \eqno{(1.11)}
$$
Then $u=0$ in $\OOO \times (0,\infty)$.
}
\\

For $0<\alpha<1$, a similar result is proved as Theorem 4.3 in 
\cite{SY}, and Theorem 2 is an improvement.

{\bf Example o $F_3$ such that $F_3\vert_{\mbox{Ker}\, (\la_n-A)}$ is 
injective.}
\\
Let  
$$
A = -\Delta, \quad d=2, \quad \OOO = \{ (x_1,x_2);\,
0<x_1<L_1,\, 0<x_2<L_2\}.         
$$
Then dim Ker $(A-\la_n) = 1$ for each $n\in \N$ if $\frac{L_1}{L_2}
\not\in \Q$.
Indeed, the eigenvalues are given by 
$\la_{mn} := \left( \frac{m^2}{L_1^2} + \frac{n^2}{L_2^2}
\right)\pi^2$, $m,n \in \N$ and the corresponding eigenfunction 
$\va_{mn}(x)$ is given by 
$\sin \frac{m\pi}{L_1}x_1 \sin \frac{n\pi}{L_2}x_2$.  
Therefore, by $\frac{L_1}{L_2} \not\in \Q$ we see that  
if $\la_{mn} = \la_{m'n'}$ with $m,n,m',n' \in \N$, then $m=m'$ and $n=n'$.

Let $x^1 = (x_1^1, x_2^1) \in \OOO$ satisfy
$\frac{x_1^1}{L_1}, \, \frac{x_2^1}{L_2} \not\in \Q$.   
We set $F_3(v) := v(x^1)$ and $M=1$.
Then we can readily verify that 
$F_3\vert_{\mbox{Ker}\, (\la_n-A)}$ is injective for all $n\in \N$.
\\

The corresponding result to Theorem 2 can be proved for the classical 
diffusion equation $\alpha=1$: 
if there exist sequences $\tau_n$, $n\in \N$ and 
$C_n>0$, $n\in \N$ which can depend on $u$ such that 
$\tau_n > 0, \quad \lim_{n\to\infty} \tau_n = \infty$ and
$$
\Vert u(\cdot,t)\Vert_{L^2(\omega)} \le C_ne^{-\tau_nt}
\quad \mbox{as $t \to \infty$},
$$
then $u=0$ in $\OOO \times (0,\infty)$.
\\

Next we consider characterizations of initial values yielding fater 
decay than $\frac{1}{t^{\alpha}}$ and/or $\frac{1}{t^{\alpha-1}}$.
\\
{\bf Theorem 3.}
\\
{\it
{\bf (i)}  Let $F_1$ be defined by (1.6).
\\
{\bf Case I: $0<\alpha<1$.}
\\
If
$$
\Vert u(\cdot,t)\Vert_{L^2(\omega)} = o\left( \frac{1}{t^{\alpha}}\right)
\quad \mbox{as $t\to \infty$},         \eqno{(1.12)}
$$
then
$$
A^{-1}a = a = 0 \quad \mbox{in $\omega$.}         \eqno{(1.13)}
$$
Moreover, assuming further that either $a\ge 0$ in $\OOO$ or 
$a\le 0$ in $\OOO$, then (1.12) yields $a=0$ in $\OOO$. 
\\
{\bf Case II: $1<\alpha<2$.}
\\
If
$$
\Vert u(\cdot,t)\Vert_{L^2(\omega)} = o\left( \frac{1}{t^{\alpha-1}}\right)
\quad \mbox{as $t\to \infty$},         \eqno{(1.14)}
$$
then
$$
A^{-1}b = b = 0 \quad \mbox{in $\omega$.}         \eqno{(1.15)}
$$
If (1.12) holds, then we have $u(x_0,0) = \ppp_tu(x_0,0) = 0$.
Moreover, assuming further that either $b\ge 0$ in $\OOO$ or 
$b\le 0$ in $\OOO$, then (1.14) yields $b=0$ in $\OOO$, and the same 
conclusion holds for $a$.
\\
{\bf (ii)}
Let $F_3$ be defined by (1.8) with $M=1$ and 
$a, b \in \DDD(A^{\gamma})$ with 
$\gamma > \frac{d}{4}$.
\\
{\bf Case 1: $0<\alpha<1$.}
\\
$$
\vert Au(x_0,t)\vert = \vert \pppa u(x_0,t)\vert 
= o\left( \frac{1}{t^{\alpha}}\right)
\quad \mbox{as $t\to \infty$}         \eqno{(1.16)}
$$
if and only if
$$
u(x_0,0) = 0.
$$
\\
{\bf Case II: $1<\alpha<2$.}
$$
\vert Au(x_0,t)\vert = \vert \pppa u(x_0,t)\vert 
= o\left( \frac{1}{t^{\alpha-1}}\right)
\quad \mbox{as $t\to \infty$}         \eqno{(1.17)}
$$
if and only if 
$$
\ppp_tu(x_0,0) = 0.
$$
Moreover (1.16) holds if and only if 
$$
u(x_0,0) = \ppp_tu(x_0,0) = 0.
$$
}
\\

Theorem 3 asserts that the faster decay than $\frac{1}{t^{\alpha}}$ or
$\frac{1}{t^{\alpha-1}}$ provides information that initial values 
vanishes at some point or in a subdomain.

In a special case, we prove
\\
{\bf Proposition 1.}
\\
{\it 
Let $a, b \in \DDD(A^{\gamma})$ with $\gamma > \frac{d}{4}$, and
$$
\left\{ \begin{array}{rl}
& a \ge 0 \quad \mbox{in $\OOO$} \quad \mbox{or} \quad 
  a \le 0 \quad \mbox{in $\OOO$}, \\
& b \ge 0 \quad \mbox{in $\OOO$} \quad \mbox{or} \quad 
  b \le 0 \quad \mbox{in $\OOO$}.
\end{array}\right.
                                                  \eqno{(1.18)}
$$
\\
{\bf Case I: $0<\alpha<1$.}
$$
\vert u(x_0,t)\vert = o\left( \frac{1}{t^{\alpha}}\right)
\quad \mbox{as $t\to \infty$},         \eqno{(1.19)}
$$
if and only if
$$
u(x,0) = 0, \quad a\in \OOO.
$$
\\
{\bf Case II: $1<\alpha<2$.}
$$
\vert u(\cdot,t)\vert = o\left( \frac{1}{t^{\alpha-1}}\right)
\quad \mbox{as $t\to \infty$}         \eqno{(1.20)}
$$
if and only if
$$
\ppp_tu(x,0) = 0, \quad x\in \OOO.
$$
}
\\

We cannot expect similar results to Theorem 2 
for the classical diffusion equation, i.e., $\alpha=1$.
\\
{\bf Example of the classical diffusion equation.}
$$
\left\{ \begin{array}{rl}
& \ppp_tu(x,t) = \ppp_x^2 u(x,t), \quad 0<x<1, \, t>0, \\
& u(0,t) = u(1,t) = 0, \quad t>0, \\
& u(x,0) = a(x), \quad 0<x<1.
\end{array}\right.
$$
Then it is well-known that for arbitrary $t_0>0$ and $a\in L^2(0,1)$,
we can choose a constant $C>0$ such that  
$$
\vert u(x_0,t)\vert \le Ce^{-\pi^2t} \quad t>0, 
$$
and
$$
\vert u(x_0,t)\vert = o(e^{-\pi^2t}) \quad \mbox{as $t\to \infty$}
$$
if and only if
$$
\sin \pi x_0 \int^1_0 a(x)\sin \pi x dx = 0.          \eqno{(1.21)}
$$
In other words, Theorem 2 means that
for $\alpha \in (0,2) \setminus \{1\}$,  
the faster decay at a point $x_0$ or in a subdomain $\omega$ still keeps 
some information of the initial value $a(x)$ at $x_0$ or in $\omega$.
On the other hand, in the case of $\alpha=1$, the decay rate is influenced 
only by averaged information (1.21) of the initial value.
However under extra assumption that the initial value $a$ does not change 
the signs, by (1.21) we can conclude that $a=0$ in $\OOO$ by 
$\sin \pi x \ge 0$ for $0<x<1$ if $\sin \pi x_0 \ne 0$.
This is true for general dimensions, because one can prove that 
the eigenfunction for $\la_1$ does not change the signs.
\\

This article is composed of five sections.  In Section 2, we show
lemmata which we use for the proofs of Theorems 1 - 3 and Proposition 1.
Sections 3 and 4 are devoted to the proofs of Theorems 1-2 and 
Theorem 3 and Proposition 1, respectivley.  In Section 5, we give concluding 
remarks.
\section{Preliminaries}

For $\alpha>0$, we define the Mittag-Leffler functions by 
$$
\MLO(z) = \sum_{k=0}^{\infty} \frac{z^k}{\Gamma(\alpha k + 1)}, \quad
\MLT(z) = \sum_{k=0}^{\infty} \frac{z^k}{\Gamma(\alpha k + 2)}, \quad z\in \C
$$
and it is know that $\MLO(z)$ and $\MLT(z)$ are entire functions in 
$z \in \C$ (e.g. Gorenflo, Kilbas, Mainardi and Rogosin \cite{GKMR},
Podlubny \cite{Po}).

First we show
\\
{\bf Lemma 1.}
\\
{\it
Let $\beta = 1,2$ and $\alpha \in (0,2) \setminus \{ 1\}$. 
\\
(i) For $p\in \N$ we have
$$
E_{\alpha,\beta}(-\eta) = \sump \frac{(-1)^{\ell+1}}
{\Gamma(\beta-\alpha\ell)}\frac{1}{\eta^{\ell}} + O\left( \frac{1}{\eta^{p+1}}
\right) \quad \mbox{as $\eta>0,\,  \to \infty$}.      \eqno{(2.1)}
$$
(ii) 
$$
\vert E_{\alpha,\beta}(\eta)\vert \le \frac{C}{1+\eta} \quad 
\mbox{for all $\eta>0$.}                                 \eqno{(2.2)}
$$
}
\\
{\bf Proof of Lemma 1.}
\\
As for (2.1), see Proposition 3.6 (pp.25-26) in \cite{GKMR} or
Theorem 1.4 (pp.33-34) in \cite{Po}.
The estimate (2.2) is seen by Theorem 1.6 (p.35) in \cite{Po} for example,
$\blacksquare$
\\

Moreover, by the eigenfunction expansion of the solution $u$ to (1.1) 
(e.g., Theorems 2.1 and 2.3 in \cite{SY}), we have
\\
{\bf Lemma 2.}
\\
{\it 
\begin{align*}
&u(x,t) = \sumn \MLO(-\la_n t^{\alpha})\sumk (a,\va_{nk})\va_{nk}(x) 
\quad \mbox{if $0 < \alpha < 1$},  \\
&u(x,t) = \sumn \biggl[ \MLO(-\la_n t^{\alpha})
\sumk (a,\va_{nk})\va_{nk}(x)
\end{align*}
$$
+ t\MLT(-\la_n t^{\alpha})
\sumk (b,\va_{nk})\va_{nk}(x) \biggr] \quad \mbox{if $1 < \alpha < 2$}
                                  \eqno{(2.4)}
$$
in $C([0,T];L^2(\OOO)) \cap C((0,T];H^2(\OOO) \cap H^1_0(\OOO))$. 
}
\\

By Lemma 1, we can prove
\\
{\bf Lemma 3.}
\\
{\it
(i) Let $a, b \in \DDD(A^{\gamma_0})$ where 
$\gamma_0=0$ if $\frac{d}{4} < 1$ and $\gamma_0 > \frac{d}{4} -1$ if 
$\frac{d}{4} \ge 1$.  Then the series in (2.4) are convergents in 
$C(\ooo{\OOO} \times [t_0,T])$.
\\
(ii) Let $a, b\in L^2(\OOO)$.  Then
\begin{align*}
&\ppp_{\nu_A}u(x,t) = \sumn \MLO(-\la_n t^{\alpha})
\sumk (a,\va_{nk})\ppp_{\nu_A}\va_{nk}(x) \quad \mbox{if $0 < \alpha < 1$},\\
&\ppp_{\nu_A}u(x,t) = \sumn \biggl[ \MLO(-\la_n t^{\alpha})
\sumk (a,\va_{nk})\ppp_{\nu_A}\va_{nk}(x)\\
+ &t\MLT(-\la_n t^{\alpha})
\sumk (b,\va_{nk})\ppp_{\nu_A}\va_{nk}(x) \biggr] \quad 
\mbox{if $1 < \alpha < 2$}
\end{align*}
in $C([t_0,T];L^2(\ppp\OOO))$. 
}
\\

For the proof of Lemma 3, we show
\\
{\bf Lemma 4.}
\\
{\it 
Let $\gamma \in \R$, and let $t_0 \in (0,T)$ 
be given arbitrarily.  We assume that $a, b \in \DDD(A^{\gamma})$.
Then there exists a constant $C=C(t_0, \gamma)>0$ such that 
$$
\Vert A^{\gamma+1}u(\cdot,t)\Vert \le
\left\{ \begin{array}{rl}
& Ct^{-\alpha}\Vert A^{\gamma}a\Vert \quad \mbox{if $0<\alpha<1$}, \\
& C(t^{-\alpha}\Vert A^{\gamma}a\Vert 
+ t^{-\alpha+1}\Vert A^{\gamma}b\Vert) \quad \mbox{if $1<\alpha<2$}
\end{array}\right.
$$
for all $t \ge t_0$.
}
\\
{\bf Proof of Lemma 4.}
\\
For $\gamma \in \R$, by each $u_0\in \DDD(A^{\gamma})$, applying (1.2) we see
\begin{align*}
& A^{\gamma+1}(u_0, \va_{nk})\va_{nk}
= (u_0, \va_{nk})\la_n^{\gamma+1}\va_{nk}
= \la_n(u_0, \la_n^{\gamma}\va_{nk})\va_{nk}
= \la_n(u_0, A^{\gamma}\va_{nk})\va_{nk}\\
=& \la_n(A^{\gamma}u_0, \va_{nk})\va_{nk}.
\end{align*}
Here we used $(u_0, A^{\gamma}\va_{nk}) = (A^{\gamma}u_0, \va_{nk})$
by (1.2).  Therefore, in view of (2.4), we have
\begin{align*}
& A^{\gamma+1}u(x,t)\\
=& \sumn \la_n \MLO(-\la_n t^{\alpha})
\sumk (A^{\gamma}a,\va_{nk})\va_{nk}(x)
+ t\sumn \la_n\MLT(-\la_n t^{\alpha}) 
\sumk (A^{\gamma}b,\va_{nk})\va_{nk}(x)
\end{align*}
in $C([0,T];L^2(\OOO))$.
We fix $t_0 > 0$ arbitrarily.
Let $1<\alpha<2$. By (2.2) we see
\begin{align*}
&\Vert A^{\gamma+1}u(\cdot,t)\Vert^2\\
\le &\sumn \la_n^2 \vert \MLO(-\la_n t^{\alpha})\vert^2
\sumk \vert (A^{\gamma}a,\va_{nk})\vert^2
+ t^2\sumn \la_n^2\vert \MLT(-\la_n t^{\alpha}) \vert^2
\sumk \vert (A^{\gamma}b,\va_{nk})\vert^2\\
\le& C(t_0)\left( \frac{1}{t^{2\alpha}}
\sumn \la_n^2 \sumk \vert (A^{\gamma}a,\va_{nk})\vert^2\frac{1}{\la_n^2}
+ \frac{1}{t^{2\alpha-2}} \sumn \la_n^2
\sumk \vert (A^{\gamma}b,\va_{nk})\vert^2\frac{1}{\la_n^2}\right)
\end{align*}
for $t \ge t_0$.  The proof for $0<\alpha<1$ is similar.
Thus we complete the proof of Lemma 4.
$\blacksquare$
\\

Now we proceed to
\\
{\bf Proof of Lemma 3 (i).}
\\
By the condition on $\gamma$, we apply the Sobolev embedding to have
$$
\Vert u(\cdot,t)\Vert_{C(\ooo{\OOO})} 
\le C\Vert A^{\gamma+1}u(\cdot,t)\Vert_{L^2(\OOO)}.
$$
Therefore, Lemma 4 yields that the series in (2.4) converge in 
$C(\ooo{\OOO} \times [t_0,T])$.  Part (ii) is seen by the trace theorem:
$$
\Vert \ppp_{\nu_A}u(\cdot,t)\Vert_{L^2(\ppp\OOO)} 
\le C\Vert Au(\cdot,t)\Vert_{L^2(\OOO)}.
$$
$\blacksquare$
\\

We conclude this section with
\\
\\
{\bf Lemma 5.}\\
{\it
We assume that $p_n \in \R$, $\{\ell_m\}_{m\in \N} \subset \N$ 
satisfying $\lim_{m\to \infty} \ell_m=\infty$, and 
there exist constants $C>0$ and $\theta_0 \ge 0$ such that 
$$
\sup_{n\in \N} \vert p_n\vert \le C\la_n^{\theta_0}.              \eqno{(2.5)}
$$
If 
$$
\sumn \frac{p_n}{\la_n^{\ell_m}} = 0 \quad \mbox{for all $m\in \N$},
$$
then $p_n=0$ for all $n\in \N$.
}
\\
{\bf Proof.}\\
By $\mu_n$, $n\in \N$, we renumber the eigenvalues $\la_n$ of $A$ according 
to the multiplicities:
$$
\mu_k = \la_1 \,\, \mbox{for $1\le k\le d_1$}, \quad
\mu_k = \la_2 \,\, \mbox{for $d_1+1\le k \le d_1+d_2$}, \cdots.
$$
Then $\mu_n\le \la_n$ for $n\in \N$.

On the other hand, 
there exists a constant $c_1>0$ such that 
$$
\mu_n = c_1n^{\frac{2}{d}} + o(1) \quad \mbox{as $n \to \infty$}
$$
(e.g., Agmon \cite{Ag}, Theorem 15.1).  Here we recall that 
$d$ is the spatial dimensions.
Therefore, we can find a constant $c_2 > 0$ such that
$\la_n \ge c_2n^{\frac{2}{d}}$ as $n \to \infty$.
Hence, we can choose a large constant $\theta_1>0$, for example
$\theta_1 > \frac{d}{2}$, such that 
$$
\sumn \frac{1}{\la_n^{\theta_1}} < \infty.
$$
We set 
$$
r_n := \frac{p_n}{\la_n^{\theta_0+\theta_1}}, \quad 
n\in \N.
$$
Then (2.5) implies 
$$
\sumn \vert r_n\vert \le \sumn \left\vert \frac{p_n}{\la_n^{\theta_0}}
\right\vert \frac{1}{\la_n^{\theta_1}}
\le C\sumn \frac{1}{\la_n^{\theta_1}} < \infty.    
$$
Since $\sumn \frac{p_n}{\la_n^{\ell_m}} = 0$, we obtain
$\sumn \frac{r_n}{\la_n^{\kappa_m}} = 0$ for all $m\in \N$, where
$\kappa_m = \ell_m - \theta_0 - \theta_1$, so that 
$$
\frac{r_1}{\la_1^{\kappa_m}} + \sum_{n=2}^{\infty} 
\frac{r_n}{\la_n^{\kappa_m}} = 0,\quad \mbox{that is},
\quad 
r_1 + \sum_{n=2}^{\infty} r_n \left( \frac{\la_1}{\la_n}
\right)^{\kappa_m} = 0.
$$
Hence 
$$
\vert r_1\vert = \left\vert -\sum_{n=2}^{\infty} r_n 
\left( \frac{\la_1}{\la_n}\right)^{\kappa_m}\right\vert 
\le \left( \sum_{n=2}^{\infty} \vert r_n\vert \right)
 \left( \frac{\la_1}{\la_2}\right)^{\kappa_m}.
$$
By $0 < \la_1 < \la_2 < ....$, we see that 
$\left\vert \frac{\la_1}{\la_2}\right\vert < 1$.
Letting $m \to \infty$, we see that $\kappa_m \to \infty$, and so
$r_1 = 0$, that is,
$p_1 =0$.  Therefore, 
$$
\sum_{n=2}^{\infty} \frac{r_n}{\la_n^{\kappa_m}} = 0.
$$
Repeating the above argument, we have $p_2=p_3= \cdots = 0$.
Thus the proof of Lemma 5 is complete.
$\blacksquare$
\section{Proofs of Theorems 1 and 2}

{\bf 3.1. Proof of Theorem 1.}

Now, by noting that $\Vert u(\cdot,t)\Vert_{H^2(\OOO)}
\le C\Vert Au(\cdot,t)\Vert$ by $u(\cdot,t) \in \DDD(A)$, 
Theorem 1 follows directly from Lemma 4 with $\gamma=0$ in Section 2.

{\bf 3.2. Proof of Theorem 2.}
\\
{\bf First Step.}
\\
It suffices to prove in the case $1<\alpha<2$, because the case 
$0<\alpha<1$ is similar and even simpler.
In view of Lemma 3, for $a$ and $b$ satisfying the conditions in the 
theorem, we have
\begin{align*}
& F_j(u(\cdot,t))
= \sumn \MLO(-\la_n t^{\alpha})
F_j\left( \sumk (a,\va_{nk})\va_{nk}\right)\\
+ & t\sumn \la_n\MLT(-\la_n t^{\alpha}) 
F_j\left( \sumk (b,\va_{nk})\va_{nk}\right), \quad j=1,2,3,4
\end{align*}
in $C([t_0,T];Y)$, where
$$
Y = 
\left\{ \begin{array}{rl}
& L^2(\omega) \qquad \mbox{for $F_1$}, \\
& L^2(\ppp\OOO) \qquad \mbox{for $F_2$}, \\
& \R^M \qquad \mbox{for $F_3$ and $F_4$}.
\end{array}\right.
$$
Applying (2.1) in Lemma 1, we obtain
$$
F_j(u(\cdot,t))
= \sump \frac{(-1)^{\ell+1}}{\Gamma(1-\alpha\ell)t^{\alpha\ell}}
\sumn \frac{1}{\la_n^{\ell}} F_j\left( \sumk (a,\, \va_{nk})\va_{nk}
\right)
+ O\left(\frac{1}{t^{\alpha p + \alpha}}\right)
\sumn F_j\left( \sumk (a,\, \va_{nk})\va_{nk}\right)
$$
$$
+ \sump \frac{(-1)^{\ell+1}}{\Gamma(2-\alpha\ell)t^{\alpha\ell-1}}
\sumn \frac{1}{\la_n^{\ell}} F_j\left( \sumk (b,\, \va_{nk})\va_{nk}\right)
+ O\left(\frac{1}{t^{\alpha p + \alpha-1}}\right)
\sumn F_j\left( \sumk (b,\, \va_{nk})\va_{nk}\right).
                                  \eqno{(3.1)}
$$
Therefore, (3.1) yields
$$
F_j(u(\cdot,t))
= \sump \frac{(-1)^{\ell+1}}{\Gamma(1-\alpha\ell)t^{\alpha\ell}}
\sumn \frac{p_n}{\la_n^{\ell}}
+ \sump \frac{(-1)^{\ell+1}}{\Gamma(2-\alpha\ell)t^{\alpha\ell-1}}
\sumn \frac{q_n}{\la_n^{\ell}} 
+ O\left(\frac{1}{t^{\alpha p + \alpha-1}}\right)\quad
                                                        \eqno{(3.2)}
$$
as $t \to \infty$.
Here we set 
$$
p_n^j = p_n = F_j\left( \sumk (a,\, \va_{nk})\va_{nk}\right), \quad
q_n^j = q_n = F_j\left( \sumk (b,\, \va_{nk})\va_{nk}\right)
$$
for $j=1,2,3,4$.

In the above series, we exclude $\ell \in \N$ such that 
$1-\alpha\ell, 2-\alpha\ell \in \{0,-1,-2, ...\}$, that is,
the terms do not appear if $\alpha \ell \in \N$.
\\
{\bf Second Step.}
\\
We see that 
$$
\mbox{$\{ \ell\in \N;\, \alpha\ell \not\in \N\}$ is an infinite set
if $\alpha \not\in \N$.}                          \eqno{(3.3)}
$$
Indeed if not, then $\left\{ \frac{n}{\alpha}\right\}_{n\in \N}
\cap \N$ is an infinite set.  Therefore there exists $N_0 \in \N$ such that 
$\left\{ \frac{n}{\alpha}\right\}_{n\in \N} \supset 
\{ N_0, N_0+1, ... \}$. Hence we can choose $n', n'' \in \N$ such that
$N_0+1 = \frac{n''}{\alpha}$ and $N_0 = \frac{n'}{\alpha}$, and so
$\frac{n''-n'}{\alpha} = 1$.  By $\alpha \not\in \N$, this is impossible.
Therefore (3.3) holds.
\\

We number the infinite set $\{ \ell \in \N;\, \alpha \ell\not\in \N\}$ by 
$\ell_1, \ell_2, \ell_3, ...$ and for each $N\in \N$, we can rewrite (3.2) as
$$
F_j(u(\cdot,t))
= \sum_{m=1}^N \frac{(-1)^{\ell_m+1}}{\Gamma(1-\alpha\ell_m)t^{\alpha\ell_m}}
\sumn \frac{p_n}{\la_n^{\ell_m}}
$$
$$
+ \sum_{m=1}^N \frac{(-1)^{\ell_m+1}}{\Gamma(2-\alpha\ell_m)t^{\alpha\ell_m-1}}
\sumn \frac{q_n}{\la_n^{\ell_m}}
+ O\left( \frac{1}{t^{\alpha\ell_{N+1}-1}}\right)
\quad \mbox{as $t\to\infty$}.                                  \eqno{(3.4)}
$$
Moreover 
$$
\{ \alpha n\}_{n \in \N} \cap \{ \alpha n-1\}_{n \in \N} 
= \emptyset \quad \mbox{for $1<\alpha<2$}.       \eqno{(3.5)}
$$
Indeed let $\alpha n' = \alpha n'' - 1$ with some $n', n'' \in \N$. Then 
$\alpha \ell_0 = 1$ with $\ell_0 := n''-n'$, which means $\alpha \le 1$ and 
this is a contradiction by $1<\alpha<2$.
\\

By (3.5), we number $\{ \alpha\ell_m\}_{m\in \N} \cup
\{ \alpha\ell_m-1 \}_{m\in \N}$ by 
$\alpha\ell_1-1=:s_1 < s_2 < \cdots < s_{2N}:= \alpha\ell_N$
and then 
$$
F_j(u(\cdot,t))
= \sum_{m=1}^{2N} \frac{Q_m}{t^{s_m}} 
+ O\left( \frac{1}{t^{\alpha\ell_{N+1}-1}}\right)
\quad \mbox{in $C([t_0,T];Y)$ as $t \to \infty$},    \eqno{(3.6)}
$$
where 
$$
Q_m = \frac{(-1)^{\ell_m+1}}{\Gamma(1-\alpha\ell_m)}
\sumn \frac{p_n}{\la_n^{\ell_m}} \quad \mbox{or} \quad
Q_m = \frac{(-1)^{\ell_m+1}}{\Gamma(2-\alpha\ell_m)}
\sumn \frac{q_n}{\la_n^{\ell_m}}.
$$
\\
{\bf Third Step.}
\\
We fix $N\in \N$ arbitrarily. 
In terms of (1.11), by (3.6) we see that for each $n\in \N$ there 
exists a constant $C_n>0$ such that 
$$
\frac{\Vert Q_1\Vert_Y}{t^{s_1}} - \sum_{m=2}^{2N} \frac{\Vert Q_m\Vert_Y}
{t^{s_m}} - \frac{C}{t^{\alpha\ell_{N+1}-1}} 
\le \frac{C_n}{t^{\tau_n}}.
$$
Then
$$
\Vert Q_1\Vert_Y \le \sum_{m=2}^{2N} \frac{\Vert Q_m\Vert_Y}
{t^{s_m-s_1}} + \frac{C}{t^{\alpha\ell_{N+1} -1-s_1}} 
+ \frac{C_n}{t^{\tau_n-s_1}}.
$$
We note that $\alpha\ell_N < \alpha\ell_{N+1} - 1$ by $\alpha>1$ and 
$\ell_n, \ell_{N+1}\in \N$, so that $s_{2N} < \alpha\ell_{N+1}-1$.

Since $\lim_{n\to\infty} \tau_n = \infty$, we can choose $n\in \N$ such that 
$\tau_n > s_1$.  Hence, letting $t \to \infty$, we have $Q_1=0$ in $Y$.
Continuing this argument, we reach $Q_m=0$ for $1 \le m \le 2N$.
Since $N\in \N$ is arbitrary, we obtain $Q_m=0$ for all $m\in \N$, that is,
$$
\sumn \frac{p_n}{\la_n^{\ell_m}} = \sumn \frac{q_n}{\la_n^{\ell_m}} = 0 
\quad \mbox{for all $m\in \N$.}
$$
In order to apply Lemma 5, we have to verify (2.5).  It suffices to 
consider for $p_n$, because the verification for $q_n$ is the same. 
\\
{\bf Case: $F_1(u(\cdot,t))$.}
\\
By the Sobolev embedding (e.g., \cite{Ad}), fixing $\mu_0>0$ with 
$2\mu_0 > d$, we have
\begin{align*}
&\Vert p_n \Vert_{C(\ooo{\OOO})}
\le C\Vert p_n\Vert_{H^{\mu_0}(\OOO)}
\le C\left\Vert A^{\frac{\mu_0}{2}} \left(
\sumk (a, \va_{nk})\va_{nk} \right) \right\Vert_{L^2(\OOO)}\\
= & C\la_n^{\frac{\mu_0}{2}} \left\Vert 
\sumk (a, \va_{nk})\va_{nk} \right\Vert_{L^2(\OOO)}
\le C\la_n^{\frac{\mu_0}{2}} \Vert a\Vert_{L^2(\OOO)}.
\end{align*}
For the second inequality, we need sufficient smoothness of the coefficients 
$a_{ij}$ and $c$ of the elliptic operator $A$ (e.g., Gilbarg and 
Trudinger \cite{GT}).  Therefore 
$$
\Vert p_n\Vert_{C(\ooo{\OOO})} \le C\la_n^{\frac{\mu_0}{2}}, \quad n\in \N.
$$
Therefore, we see (2.5) for $F_1, F_3$ and $F_4$ with 
$\theta_0 = \frac{\mu_0}{2}$.
\\
{\bf Case: $F_2(u(\cdot,t))$.}
\\
We fix $\mu_0>0$ such that $2\mu_0 > d$. Then by the Sobolev embedding,
we obtain
\begin{align*}
& \left\Vert \ppp_{\nu_A}\left(\sumk (a, \va_{nk})\va_{nk} \right)
\right\Vert_{C(\ppp\OOO)}
\le C\left\Vert \left(
\sumk (a, \va_{nk})\va_{nk} \right) \right\Vert_{C^1(\ooo{\OOO})}
\le C\left\Vert 
\sumk (a, \va_{nk})\va_{nk} \right\Vert_{H^{\mu_0+1}(\OOO)}\\
\le& C\left\Vert A^{\frac{\mu_0}{2}+\frac{1}{2}} 
\sumk (a, \va_{nk})\va_{nk} \right\Vert_{L^2(\OOO)}
= C\la_n^{\frac{\mu_0}{2}+\frac{1}{2}} \left\Vert 
\sumk (a, \va_{nk})\va_{nk} \right\Vert_{L^2(\OOO)}
\le C\la_n^{\frac{\mu_0}{2}+\frac{1}{2}} \Vert a\Vert_{L^2(\OOO)}.
\end{align*}
Hence (2.5) is satisfied with $\theta_0 = \frac{\mu_0}{2} + \frac{1}{2}$.
\\

Therefore, Lemma 5 yields $p_n=q_n = 0$ for all $n\in \N$, that is,
$$
F_j\left( \sumk (a,\, \va_{nk})\va_{nk} \right)
= F_j\left( \sumk (b,\, \va_{nk})\va_{nk} \right) = 0,
\quad j=1,2,3,4, \quad n\in \N.                        \eqno{(3.7)}
$$
\\
{\bf Fourth Step.}
\\
It suffices to verify that $p_n=0$ for $n\in \N$ imply $a=0$ in $\OOO$.
For $F_3$ and $F_4$, the assumption in Theorem 2 yields 
$$
\sumk (a, \,\va_{nk})\va_{nk} = \sumk (b, \,\va_{nk})\va_{nk} = 0
\quad \mbox{in $\OOO$}
$$
for all $n\in \N$.  Therefore, $a=b=0$ in $\OOO$, that is,
$u=0$ in $\OOO\times (0,\infty)$.  Thus the proof of Theorem 2 is 
complete for $F_3$ and $F_4$.
\\
{\bf Case: $F_1$.}
By (3.7), we have
$$
p_n(x) = \sumk (a,\, \va_{nk})\va_{nk}(x) = 0, \quad
n\in \N, \, x\in \omega.
$$
Since $(A-\la_n)p_n = 0$ in $\OOO$, we apply the unique
continuation for the elliptic operator $A-\la_n$ (e.g., Choulli 
\cite{Ch}, H\"ormander \cite{Ho}) to see that $p_n=0$ in $\OOO$ for $n\in \N$.  Since $a= \sumn p_n$ in $L^2(\OOO)$, we reach
$a = 0$ in $\OOO$.
\\
{\bf Case: $F_2$.}
We set $u_n(x) = \sumk (a,\va_{nk})\va_{nk}(x)$ for $x\in \OOO$.
By $u_n \in \mathcal{D}(A)$, we have $u_n=0$ on $\Gamma$ and so
$$
\ppp_{\nu_A}u_n(x) = u_n(x) = 0, \quad n\in \N, \, x\in \Gamma.
$$
Therefore, since $(A-\la_n)u_n = 0$ in $\OOO$, the unique continuation 
(e.g.,  \cite{Ch}, \cite{Ho}) yields\\
$\sumk (a,\va_{nk})\va_{nk}(x) = 0$ for all $n\in \N$ and $x\in \Omega$.
Hence, we can see $a=0$ in $\OOO$.  
Thus the proof of Theorem 2 is complete.
\section{Proofs of Theorem 3 and Proposition 1}

{\bf 4.1. Proof of Theorem 3}
\\
{\bf Case: $F_1$.}
\\
It is sufficient to prove the case $1<\alpha<2$.
Let (1.14) hold.  By (3.2) with $p=1$, noting that $\Gamma(1-\alpha)$
and $\Gamma(2-\alpha)$ are finite, we see
$$
\left\Vert \frac{1}{\Gamma(1-\alpha)}\frac{1}{t^{\alpha}}
\sumn \sumk \frac{(a,\,\va_{nk})\va_{nk}}{\la_n}
+ \frac{1}{\Gamma(2-\alpha)}\frac{1}{t^{\alpha-1}}
\sumn \sumk \frac{(b,\,\va_{nk})\va_{nk}}{\la_n}\right\Vert_{L^2(\omega)}
= o\left( \frac{1}{t^{\alpha-1}}\right).
                                                     \eqno{(4.1)}
$$
Therefore, in terms of (1.3), we obtain
$$
\frac{1}{\Gamma(2-\alpha)}\frac{1}{t^{\alpha-1}}\Vert A^{-1}b\Vert
_{L^2(\omega)}
- \frac{1}{\Gamma(1-\alpha)}\frac{1}{t^{\alpha}}\Vert A^{-1}a\Vert
_{L^2(\omega)}
= o\left( \frac{1}{t^{\alpha-1}}\right)
$$
as $t\to \infty$.
Multiplying with $t^{\alpha-1}$ and letting $t\to\infty$, we obtain
$A^{-1}b=0$ in $\omega$.  

Next let (1.12) hold.  Then, by $o\left( \frac{1}{t^{\alpha}}\right)
\le o\left( \frac{1}{t^{\alpha-1}}\right)$, we have also (1.14), so that
we have already proved $A^{-1}b=0$ in $\omega$.
Therefore, since 
$$
\frac{(-1)^2}{\Gamma(2-\alpha)}\frac{1}{t^{\alpha-1}}
\sumn \sumk \frac{(b,\,\va_{nk})\va_{nk}}{\la_n} 
= \frac{1}{\Gamma(2-\alpha)}\frac{1}{t^{\alpha-1}} A^{-1}b = 0 \quad 
\mbox{in $\omega$},
$$
equality (3.2) with $p=1$ and (1.12) yield 
$$
\frac{1}{\Gamma(1-\alpha)}\frac{1}{t^{\alpha}}\Vert A^{-1}a\Vert_{L^2(\omega)}
+ o\left( \frac{1}{t^{2\alpha}}\right)
+ o\left( \frac{1}{t^{2\alpha-1}}\right)
= o\left( \frac{1}{t^{\alpha}}\right).
$$
Multiplying with $t^{\alpha}$ and letting $t\to \infty$, by $\alpha-1>0$,
we see that $A^{-1}a=0$ in $\omega$.

Moreover $A^{-1}a=0$ in $\omega$ implies $a=0$ in $\omega$.  Indeed,
setting $g:= A^{-1}a$ in $\OOO$, we have $g=0$ in $\omega$ and 
$Ag = a$ in $\OOO$. Therefore, $a=A0=0$ in $\omega$.  Similarly
$A^{-1}b=0$ in $\omega$ yields $b=0$ in $\omega$.
\\

Finally we have to prove that the extra condition
$$
\mbox{$a\ge 0$ in $\OOO$ or $a\le 0$ in $\OOO$,}    \eqno{(4.2)}
$$
implies $a=0$ in $\OOO$.

Let $a\ge 0$ in $\OOO$.  Then $g:= A^{-1}a$ satisfies 
$$
\sum_{i,j=1}^d \ppp_i(a_{ij}(x)\ppp_jg(x)) + c(x)g(x) \ge 0 
\quad \mbox{in $\OOO$}.
$$
By $c\le 0$ in $\OOO$ and $g=0$ on $\ppp\OOO$, the weak maximum principle 
(e.g., Theorem 3.1 (p.32) in Glilbarg and Trudinger \cite{GT})
implies that $g \le 0$ on $\ooo{\OOO}$.
Since $g(x) = 0$ for $x\in \omega$, we see that $g$ achieves 
the maximum $0$ at an interior point $x_0\in \OOO$.  
Again by $c\le 0$ in $\OOO$, the strong maximum principle (e.g., 
Theorem 3.5 (p.35) in \cite{GT}) yields that 
$g$ is a constant function, that is, $g(x) = 0$ for all $x \in \OOO$.
Hence, $a=Ag=0$ in $\OOO$.
Thus the proof in the case $F_1$ is complete.
\\
{\bf Case: $F_3$.}
\\
It suffices to prove only in the case $1<\alpha<2$.  By Lemma 2,
for arbitrarily chosen $t_0 \in (0,T)$, we see 
\begin{align*}
& Au(x,t) 
= \sumn E_{\alpha,1}(-\la_nt^{\alpha}) \sumk (a,\, \va_{nk})\la_n\va_{nk}\\
+ &\sumn tE_{\alpha,2}(-\la_nt^{\alpha}) \sumk (b,\, \va_{nk})\la_n\va_{nk}
\quad \mbox{in $C([t_0,T];L^2(\OOO))$}.
\end{align*}
Using $a, b\in \DDD(A^{\gamma})$ with $\gamma > \frac{d}{4}$ and 
$$
A^{\gamma}(a, \,\va_{nk})\la_n\va_{nk} = \la_n^{1+\gamma}(a, \,\va_{nk})
\va_{nk}
= \la_n(a, \,A^{\gamma}\va_{nk})\va_{nk}
= \la_n(A^{\gamma}a, \, \va_{nk})\va_{nk}, \quad \mbox{etc.,}
$$
we obtain
\begin{align*}
& A^{1+\gamma}u(x,t) 
= \sumn E_{\alpha,1}(-\la_nt^{\alpha}) \sumk (A^{\gamma}a,\, \va_{nk})
\la_n\va_{nk}\\
+ &\sumn tE_{\alpha,2}(-\la_nt^{\alpha}) \sumk (A^{\gamma}b,\, \va_{nk})
\la_n\va_{nk}.
\end{align*}
Consequently, by Lemma 1, we can prove  
$$
\Vert A^{1+\gamma}u\Vert_{L^{\infty}(t_0,T;L^2(\OOO))} < \infty,
$$
and so the above series is convergent in $L^{\infty}(t_0,T;L^2(\OOO))$.
Since the Sobolev embedding implies $\DDD(A^{\gamma}) \subset 
C(\ooo{\OOO})$ with $\gamma > \frac{d}{4}$, we obtain
\begin{align*}
& Au(x_0,t) 
= \sumn E_{\alpha,1}(-\la_nt^{\alpha}) \sumk (a,\, \va_{nk})
\la_n\va_{nk}(x_0)\\
+ &\sumn tE_{\alpha,2}(-\la_nt^{\alpha}) \sumk (b,\, \va_{nk})
\la_n\va_{nk}(x_0), \quad t_0< t < T \quad \mbox{in $C[t_0,T]$}.
\end{align*}
Substituting (2.1) with $p=1$ and $\beta = 1,2$, we have
\begin{align*}
& Au(x_0,t) 
= \frac{1}{\Gamma(1-\alpha)} \sumn \sumk (a,\, \va_{nk})\va_{nk}(x_0)
\frac{1}{t^{\alpha}}\\
+& \frac{1}{\Gamma(2-\alpha)} \sumn \sumk (b,\, \va_{nk})\va_{nk}(x_0)
\frac{1}{t^{\alpha-1}}
+ O\left( \frac{1}{t^{2\alpha-1}}\right) \quad 
\mbox{as $t \to \infty$}.
\end{align*}
By $a, b \in \DDD(A^{\gamma}) \subset C(\ooo{\OOO})$, we find
$$
Au(x_0,t) = \frac{1}{\Gamma(1-\alpha)t^{\alpha}} a(x_0)
+ \frac{1}{\Gamma(2-\alpha)t^{\alpha-1}} b(x_0)
+ O\left( \frac{1}{t^{2\alpha-1}}\right) \quad \mbox{as $t \to \infty$}.
                               \eqno{(4.3)}
$$
By an argument similar to Case $F_1$ in Theorem 3, we see that 
(1.16) and (1.17) imply $a(x_0) = 0$ and $b(x_0) = 0$ respectively.
The converse assertion in the theorem directly follows from (4.3).

{\bf 4.2. Proof of Proposition 1}

It is sufficient to prove in the case $1<\alpha<2$. 
By $a, b \in \DDD(A^{\gamma}) \subset C(\ooo{\OOO})$ with 
$\gamma > \frac{d}{4}$, similarly to (4.1), we obtain
\begin{align*}
& u(x_0,t) \\
= & \frac{1}{\Gamma(1-\alpha)}\frac{1}{t^{\alpha}}
\sumn \sumk \frac{(a,\, \va_{nk})}{\la_n}\va_{nk}(x_0)
+ \frac{1}{\Gamma(2-\alpha)}\frac{1}{t^{\alpha-1}}
\sumn \sumk \frac{(b,\, \va_{nk})}{\la_n}\va_{nk}(x_0)
+ O\left( \frac{1}{t^{2\alpha-1}}\right)\\
= &\frac{1}{\Gamma(1-\alpha)}\frac{1}{t^{\alpha}}(A^{-1}a)(x_0)
+ \frac{1}{\Gamma(2-\alpha)}\frac{1}{t^{\alpha-1}}(A^{-1}b)(x_0)
+ O\left( \frac{1}{t^{2\alpha-1}}\right) \quad \mbox{as $t \to \infty$}.
\end{align*}
Similarly to the case $F_1$ in the proof of Theorem 3, we can prove 
that (1.19) and (1.20) imply $(A^{-1}a)(x_0) = 0$ and
$(A^{-1}b)(x_0) = 0$ respectively.
Under the assumption that $a$ and $b$ do not change the signs in $\OOO$,
in view of the weak and the strong maximum principles, we can 
argue similarly to the final part of the proof of Theorem 3 in the 
case of $F_1$, so that we can reach $a=0$ in $\OOO$ and/or $b=0$ in $\OOO$.
Therefore, we prove that (1.19) and (1.20) imply $a(x) = b(x) = 0$ 
and $b(x) = 0$ for $x\in \OOO$, respectively.
The converse statement of the proposition is readily seen.
Thus the proof of Proposition 1 is complete.

\section{Concluding remarks}

{\bf 5.1.}
\\
Time-fractional diffusion-wave equations with order
$\alpha\in (0,2) \setminus \{1 \}$ describe slow diffusion and is known not 
to have strong smoothing property as the classical diffusion equation.
Such a weak smoothing property is characterized by the norm 
equivalence between $\Vert u(\cdot, t)\Vert_{H^2(\OOO)}$ and
$\Vert u(\cdot,0)\Vert_{L^2(\OOO)}$ for any $t>0$ in the case of 
$0<\alpha<1$. 
The weak smoothing property allows that the backward problem in time 
is well-posed for $\alpha \in (0,2) \setminus \{1\}$ 
(Floridia, Li and Yamamoto \cite{FLY}, Floridia and Yamamoto \cite{FY},
Sakamoto and Yamamoto \cite{SY}), which is a remarkable difference from 
the case $\alpha=1$.

The current article establishes that local properties of initial values
affect the decay rate of solution as $t \to \infty$, which indicates 
that a time-fractional equation can keep some profile of the initial 
value even for very large $t>0$, which can be understood related to the 
backward well-posedness in time and is essentially different from 
the case $\alpha=1$.

The essence of the argument relies on that the behavior 
of a solution $u$ for large $t>0$ admits an asymptotic 
expansion with respect to $\left( \frac{1}{t}\right)^{\alpha\ell}$
and $\left( \frac{1}{t}\right)^{\alpha\ell-1}$ with $\ell \in \N$.
\\

{\bf 5.2.}
\\
We can generalize Theorem 3 (ii).  For simplicity, we consider only 
the case $0<\alpha<1$.
\\
{\bf Proposition 2.}
\\
{\it 
Let $a \in \DDD(A^{\gamma})$ with $\gamma > \frac{d}{4}$ and 
$0 < \alpha, \beta < 1$.  Then
\\
(i) 
$$
\vert \ppp_t^{\beta}u(x_0,t)\vert \le \frac{C}{t^{\beta}}\Vert a\Vert.
$$
(ii) If 
$$
\vert \ppp_t^{\beta}u(x_0,t)\vert = o\left(\frac{1}{t^{\beta}}\right)
\quad \mbox{as $t \to \infty$}, 
$$
then $u(x_0,0) = 0$.
}
\\
The proof relies on 
$$
\ppp_t^{\beta}u(x,t) 
= -t^{\alpha-\beta}\sumn \la_n E_{\alpha,\alpha+1-\beta}
(-\la_nt^{\alpha}) \sumk (a, \,\va_{nk)}\va_{nk}(x)
\quad \mbox{in $C((0,T];L^2(\OOO))$}                     \eqno{(5.1)}
$$
and then we can argue similarly to Theorem 3 (ii) by (2.1).
The equation (5.1) can be verified as follows:
$$
\ppp_t^{\beta}(t^{\alpha k}) = \frac{\Gamma(\alpha k + 1)}
{\Gamma(\alpha k + 1 - \beta)}t^{\alpha k - \beta}, \quad 
k\in \N,
$$
and so the termwise differentiation yields
$$
\ppp_t^{\beta}E_{\alpha,1}(-\la_nt^{\alpha})
= -\la_n t^{\alpha-\beta} E_{\alpha,\alpha+1-\beta}(-\la_nt^{\alpha}),
\quad t>0.
$$
Then (2.4) yields (5.1).

We omit the details of the proof of Proposition 2.
\section*{Acknowledgment}
The author was supported by Grant-in-Aid for Scientific Research (S)
15H05740 and Grant-in-Aid (A) 20H00117 of 
Japan Society for the Promotion of Science and
by The National Natural Science Foundation of China
(no. 11771270, 91730303).
This paper has been supported by the RUDN University 
Strategic Academic Leadership Program.


\end{document}